\documentclass[12pt]{article}
\usepackage{graphics,latexsym,amssymb,amsmath,amscd}
\RequirePackage{graphics}
\def\[#1\]{\begin{eqnarray*}#1\end{eqnarray*}}
\def\Re{\hbox{\rm Re}\,}

\def\phi{\varphi}

\input xy
\xyoption{all}
\newtheorem{thm}{Theorem}[section]
\newtheorem{dfn}[thm]{Definition}

\newtheorem{prop}[thm]{Proposition}
\newtheorem{lemma}[thm]{Lemma}
\newtheorem{ex}[thm]{Example}

\newcommand{\Pf}{{\em Proof}. }
\newcommand{\EPf}{\hbox{}\hfill$\Box$\vspace{.5cm}}
\newcommand{\C}{{{\mathbb C}}}
\newcommand{\R}{{{\mathbb R}}}

\def\Re{\hbox{\rm Re}\,}

\def\phi{\varphi}

\def\C{{\mathbb C}}
\def\R{{\mathbb R}}

\date{}
\title{Diffeomorphisms preserving $\R$-circles in three dimensional CR manifolds}
\author{E. Falbel \\Institut de Math{\'e}matiques\\
 Analyse Alg\'ebrique case 82\\ F-75252 Paris France\\ \and
 J. M. Veloso\\  Faculdade de Matem\' atica - ICEN\\
Universidade Federal do Par\'a\\66059 - Bel\' em- PA - Brazil}
\begin{document}
\maketitle
\newcommand{\D}{\mbox{$\cal D$}}
\newtheorem{df}{Definition}[section]
\newtheorem{te}{Theorem}[section]
\newtheorem{co}{Corollary}[section]
\newtheorem{po}{Proposition}[section]
\newtheorem{lem}{Lemma}[section]
\newcommand{\Ad}{\mbox{Ad}}
\newcommand{\ad}{\mbox{ad}}
\newcommand{\im}[1]{\mbox{\rm im\,$#1$}}
\newcommand{\bm}[1]{\mbox{\boldmath $#1$}}
\newcommand{\sime}{\mbox{sim}}
\begin{abstract}
$\R$-circles in general three dimensional CR manifolds (of contact type)  are the analogues to traces of Lagrangian totally geodesic planes on $S^3$ viewed as 
the boundary of  two dimensional complex hyperbolic space.  They form a family of certain legendrian curves on the manifold.
We prove that a diffeomorphism between three dimensional CR manifolds which preserve circles is
either a CR diffeomorphism or conjugate CR diffeomorphism. 
\end{abstract}
\section{Introduction}

Given a  three manifold $M$ equipped with a contact plane distribution $D$, we say $M$ is a CR manifold if $D$ is equipped with 
a complex operator $J : D\rightarrow D$ satisfying $J^2=-Id$.   Examples of CR manifolds arise naturally as real hypersurfaces in $\C^2$ or as boundaries 
of complex two-dimensional
manifolds.  The most important example appears as follows:
the complex two-dimensional ball 
$$
H^2_\C=\left \{ \ (z_1,z_2)\in \C^2\ \vert \ \vert z_1\vert^2+ \vert z_2\vert^2 < 1\ \right\}
$$
has boundary $S^3$ which is naturally equipped with a distribution $D=TS^3\cap JTS^3$, where $J$ denotes the standard complex structure of $\C^2$.
The group of biholomorphisms of $H^2_\C$ is $PU(2,1)$ and it acts on $S^3$.

$\C$-circles (or chains) and  $\R$-circles (or pseudocircles or circles for short) are analogues defined for general CR manifolds to curves in $S^3$ obtained as traces
 of complex lines and Lagrangian totally geodesic subspaces of
 ${H^2_\C}$ in its boundary. 
  More precisely, complex lines through the origin are totally geodesic spaces and
their intersection with $S^3$ are called chains.  The group of biholomorphisms acts transitively on the space of chains.  It turns out that a Lagrangian plane passing through the origin is a totally
real totally geodesic surface of  ${H^2_\C}$.  Its intersection with the sphere is called an $\R$-circle.  Again, $PU(2,1)$ acts transitively on the space of 
Lagrangian totally geodesic subspaces.

The generalization of these curves to  real hypersurfaces of $\C^2$ was considered by E. Cartan \cite{C}, where chains and circles (what we call here $\R$-circles) are defined.   Here we will restrict our attention to $\R$-circles.
Cartan's  definition was considered by Jacobowitz in \cite{J}, chapter 9, where these curves are called pseudocircles and are 
 defined on an abstract three dimensional CR manifold. 
In this paper we make explicit the system of differential equations satisfied by $\R$-circles, and use these equations to prove our main theorem (Theorem \ref{main}) that a diffeomorphism between two CR manifolds which preserve circles is a CR map or an anti CR map. Our motivation was the analogue result for chains in \cite{Ch}.

We will first recall general results on CR structures, the group $SU(2,1)$ and the Cartan's connection on a principal bundle $Y\rightarrow M$ with
group $H$ (a subgroup of $SU(2,1)$ which is the isotropy of its action on $S^3$) which will be used in our main theorem.  We also recall
the definition of $\R$-circles (also called pseudocircles or circles for short) as projections of leaves of an integrable differential system (for more details we refer to \cite{J}).  We finally obtain the differential equations
which circles satisfy in $M$ (Theorem \ref{eqcircle}) and use them in the proof of our main theorem in the last section.

\section{CR structures}

For a general reference or this section see \cite{CM,J}.
We consider a three manifold $M$ equipped with a contact plane distribution $D$.
A 3-dimensional \emph{CR-structure} $(M,D,J)$ is the contact manifold $M$ equipped with the complex operator $J : D\rightarrow D$ satisfying $J^2=-Id$. 

If $(M,D,J)$ and $(\tilde M,\tilde D,\tilde J)$ are 3-dimensional CR-structures, then a diffeomorphism 
$f:M\rightarrow \tilde M$ is a \emph{CR-diffeomorphism} if $f_*(D)=\tilde D$ and $f_*J=\tilde Jf_*$. If $f_*J=-\tilde Jf_*$, we say that $f$ is a \emph{conjugate}
 CR-diffeomorphism.

A CR-structure induces an orientation
on $D$ and an orientation of the normal bundle $TM/TD$ given by
$X,JX,[X,JX]$ where $X$ is a local section of $D$. 

Fixing a local section $X$ of $D$ one can define a form $\theta$ such that $\theta(D)=0$ and such that
$\theta([X,JX])=-2$.

Consider $D\otimes \C=D^{1,0}\oplus D^{0,1}\subset TM\otimes \C$.  
Taking $Z=\frac{1}{2}(X-iJX)$ and $\bar Z=\frac{1}{2}(X+iJX)$, we get $[Z,\bar Z]=\frac{i}{2}[X,JX]$. We define a form  $\theta^1\in {D^{0,1}}^\perp$ such that $\theta^1(Z)=1$. Then 
$d\theta(Z,\bar Z)=-\theta([Z,\bar Z])=i$, so
$$
d\theta =i\theta^1\wedge \theta^{\bar 1} \ {\mbox {modulo }}\theta
$$
where we define $\theta^{\bar 1}=\overline{\theta^{ 1}}$.  If $\theta'^1$
is another form satisfying the equation we have
$$
\theta'^1=e^{i\alpha}\theta^1 \ {\mbox {modulo }}\theta
$$
for $\alpha\in \R$.

 Let ${E}$ to be the oriented line bundle of all
forms $\theta$ as above.  On $E$ we define the tautological
form $\omega$.  That is $\omega_\theta=\pi^*(\theta)$ where
$\pi: {E}\rightarrow M$ is the natural projection.

We consider the tautological forms defined
by the forms above over the line bundle $E$.  That is, for each 
$\theta^1$ as above, we let $\omega^1_\theta= \pi^*(\theta^1)$.
At each point $\theta\in E$ we have the family
$$
        \omega' = \omega
$$
$$
        \omega'^{1} = e^{i\alpha} \omega^1 + v^1\omega
$$
where we understand that the forms are defined over $E$.
  Those forms
vanish on vertical vectors, that is, vectors in the kernel
of the map $TE\rightarrow TM$. In order to define non-horizontal 1-forms 
we let $\theta$ be a section of $E$ over $M$ and introduce the coordinate
$\lambda\in \R^+$ in $E$.
By abuse of notation,
let $\theta$ denote the tautological form on the section $\theta$.
Therefore the tautological form $\omega$ over $E$ is 
$$
\omega_{\lambda}=\lambda \theta.
$$
Differentiating this formula we obtain
\begin{equation}
                d\omega = \omega\wedge\phi + 
i \omega^1\wedge\omega^{\bar 1}  \label{domega}
\end{equation}
where $\phi= -\frac{d\lambda}{\lambda}$ 
 modulo $\omega, \omega^1, 
\omega^{\bar 1}$, $\phi$ real.

 Observe that
$
\frac{d\lambda}{\lambda}
$ is a form intrinsically defined on $E$ up to horizontal
forms (the minus sign is just a matter of conventions).
 In fact, choosing a different section $\theta'$ with 
$\theta=\mu \theta'$ where  $\mu$ is a function over $M$, we can write
$\omega={\lambda}\theta={\lambda\mu}\theta'$ and obtain
$$
\frac{d(\lambda\mu)}{\lambda\mu}=
\frac{\mu d\lambda+\lambda d\mu}{\lambda\mu}=
\frac{d\lambda}{\lambda}+\frac{d\mu}{\mu},
$$
where $\frac{d\mu}{\mu}$ is a horizontal form.

If we write (\ref{domega}) as 
$$\begin{array}{rcl}
d\omega&=&\omega'\wedge\phi'+i\omega'^{1}\wedge \omega'^{\bar 1}=\omega\wedge \phi'+i(e^{i\alpha} \omega^1 + v^1\omega)\wedge (e^{-i\alpha} \omega^{\bar 1} + v^{\bar 1}\omega)\\
&=&
\omega\wedge (\phi'-ie^{i\alpha}v^{\bar 1}\omega^1 +ie^{-i\alpha}v^1\omega^{\bar 1})+i \omega^1\wedge\omega^{\bar 1}
\end{array}
$$
it follows that  $\phi'-ie^{i\alpha}v^{\bar 1}\omega^1 +ie^{-i\alpha}v^1\omega^{\bar 1}=\phi  \ {\mbox {modulo }}\omega$.

We obtain in this way a coframe bundle  over $E$:
$$
        \omega' = \omega
$$
$$
        \omega'^{1} = e^{i\alpha} \omega^{1} + v^1 \omega
$$
$$
        \phi'= \phi + ie^{i\alpha}v^{\bar 1}\omega^1- 
ie^{-i\alpha}v^{ 1}\omega^{\bar 1} +s\omega
$$
$v^1\in \C$  and $s\in \R$ are arbitrary.

\begin{dfn}
We denote by $Y$ the  coframe bundle
 $Y\rightarrow E$ given by the set of 1-forms  $\phi, \omega^{1},
  \omega^{\bar 1}, \omega$.
  Two coframes
are related by
$$
(\phi', \omega'^{1},
  \omega'^{\bar 1}, \omega')=
(\phi, \omega^{1},
  \omega^{\bar 1}, \omega)
 \left ( \begin{array}{cccc}

                        1       &    0                  &       0 & 0 \\

                         ie^{i\alpha}v^{\bar 1}   &     e^{i\alpha} & 0            &       0 \\
                        -ie^{-i\alpha}v^{ 1}   &  0 & e^{-i\alpha}             &       0 \\
                        s       &  v^1 & v^{\bar 1}     &       1

                \end{array} \right )
$$
where   $s, \alpha \in \R $ and $v^1\in \C$.
\end{dfn}

\begin{thm}[\cite{C,CM}]\label{structure}
On $Y$ we have unique globally defined forms $\omega,\omega^1,\omega^1_1,\phi^1,\psi$ such that
$$
\begin{array}{rcl}
 d\omega &= &\omega\wedge\phi + 
i \omega^1\wedge\omega^{\bar 1} \vspace{.3cm}\\
d\omega^1&=&\frac{1}{2} \omega^1\wedge \phi +
\omega^1\wedge\omega^1_1 +\omega\wedge \phi^1\vspace{.3cm}\\
d\phi &= &
i\omega ^1\wedge\phi^{\bar 1} -
i\omega^{\bar 1}\wedge \phi^1 +\omega \wedge \psi\vspace{.3cm}\\
d\omega^1_1&=&\frac{3}{2}i\omega^1\wedge \phi^{\bar 1}+
\frac{3}{2}i\omega^{\bar 1}\wedge \phi^{ 1}\vspace{.3cm}\\
d\phi^1&=&\frac{1}{2}\phi\wedge \phi^1-\omega^1_1\wedge\phi^1
+\frac{1}{2}\omega^1\wedge \psi+ Q^1_{\bar 1}\omega\wedge \omega^{\bar 1}\vspace{.3cm}\\
d\psi&=& 2i\phi^1\wedge \phi^{\bar 1}+\phi\wedge\psi+
\left ( R_{ 1}\omega^{ 1}+R_{\bar 1}\omega^{\bar 1} \right)\wedge \omega
\end{array}
$$
with $\omega^1_1+\overline{\omega^1_1}=0$, $\omega, \phi,\psi$ real and
$$
dQ^1_{\bar 1}-2Q^1_{\bar 1}\phi
+2Q^1_{\bar 1}\omega^1_1 = 
S\omega -\frac{1}{2}R_{\bar 1}\omega^{ 1}+T\omega^{\bar 1},
$$
$$
dR_1-\frac{5}{2}R_1\phi-R_1\omega^1_1+
2i\overline Q^1_{\bar 1}\phi^1 =A\omega+B\omega^1+C\omega^{\bar 1}.
$$
where $A,B,C$ are functions on $Y$ and $C$ is real.
\end{thm}

We can verify easily that structure equations in Cartan or Jacobowitz are the same as here, with the correspondance
$$\Omega=\omega,\ \,\Omega_1=\omega^1,\ \,\Omega_2=-\frac{1}{2}\phi-\omega^1_1,\ \,\Omega_3=-\phi^1,\ \,\Omega_4=\frac{1}{2}\psi$$
and
$$R=Q^1_{\bar 1},\ \,S=R_1.$$
\section{${\bf SU(2,1)}$}

Define $$ SU(2,1)= \{\ g\in SL(3,\C)\ |\ \bar{g}^T Q g=Q\ \} $$
where the Hermitian form $Q$  given by 
$$
 Q =\left ( \begin{array}{ccc}
                             0  & 0 & i/2\\
                             0  & 1 & 0\\
                             -i/2  & 0 & 0

               \end{array}
       \right ).
$$ 
The group $SU(2,1)$ acts on $\C^{3}$ on the left preserving
the cone 
$$ \{\ z\in \C^{3}\ |\ \bar{z}^T Q z=0\ \}.$$
  The
projectivization of this cone is $S^{3}\subset \C P^{2}$.
$SU(2,1)$ 
has a finite  center
$K$ which is a cyclic group of order $3$ 
acting trivially on the sphere  $S^{3}$. We
define $PU(2,1)= SU(2,1)/K$.

The elements of the Lie algebra
$ {\mathfrak su(2,1)}$ are represented by the matrices $$
   \left ( \begin{array}{ccc}
                        u      &    -2i\bar{y}    &   w   \\

                        x     &   a     &  y  \\    

                        z       &  2i\bar{x}   &       -\bar{u}
                \end{array} \right )
$$ where $ia\in \R$, $z,w\in \R$, $x,y\in \C$,
$u\in\C$ and  $u-\bar{u}=-a$.
Observe that the Lie algebra $\mathfrak g=\mathfrak su(2,1)$ is graded:
$$
\mathfrak  g
= \mathfrak g^{-2}\oplus \mathfrak g^{-1}\oplus \mathfrak g^0\oplus \mathfrak g^1\oplus \mathfrak g^2 $$
where 
$$
\mathfrak g^{-2}= \left\{\left ( \begin{array}{ccc}

                        0      &    0    &   0   \\

                        0     &    0     & 0 \\

                        z       &  0    &     0
                \end{array} \right )\right\}\ \ \
 \mathfrak g^{-1}= \left\{\left ( \begin{array}{ccc}

                        0      &    0    &   0   \\

                        x     &    0     &   0 \\

                        0       &  2i\bar{x}    &     0

                \end{array} \right )\right\}
$$ 
$$ 
\mathfrak g^{0}=\left\{\left ( \begin{array}{ccc}

                        u      &    0    &  0   \\

                        0    &    a     &   0\\

                        0       &  0    &       -\bar{u}

                \end{array} \right )\right\}\ \ \
\mathfrak g^{1}=\left\{\left ( \begin{array}{ccc}

                        0      &     -2i\bar{y}    &   0   \\
                       0     &    0     &   y \\

                        0       &  0    &     0

                \end{array} \right )\right\}\ \ \
\mathfrak g^{2}=\left\{\left ( \begin{array}{ccc}

                        0      &    0    &   w   \\

                        0   &    0     &   0 \\

                        0       &  0    &     0

                \end{array} \right )\right\}
$$
We have
 $$ \mathfrak g^{0}=\R\oplus \mathfrak u(1),$$
where $$ \mathfrak u(1)=\left\{\left (
\begin{array}{ccc}

                        -iq/2      &    0    &  0   \\

                        0    &    iq     &   0\\

                        0       &  0    &       -iq/2

                \end{array} \right )\right\}
$$with  $q\in\R$ and $$\R=\left\{\left (
\begin{array}{ccc}

                        r      &    0    &  0   \\

                        0    &    0     &   0\\

                        0       &  0    &       -r

                \end{array} \right )\right\} $$ with $r\in\R$.

Define the subalgebra $$ \mathfrak h= \mathfrak g^0 \oplus \mathfrak g^1\oplus \mathfrak g^{2}.$$

The isotropy of the action of $SU(2,1)/K$ on $S^{3}$ at the
point $[1,0 , 0]^{T}$ is the group $H=CU(1) \ltimes N$
(whose Lie algebra is $\mathfrak h$),  represented (up to $K$) by matrices of the form 
\begin{equation}\label{h}
   \left ( \begin{array}{ccc}

                        a       &    -2i\bar{a}\bar{b}    &       a(s-ib\bar{b})  \\

                        0   &     \frac{\bar{a}}{a}    &       b \\

                        0      &  0    &       \bar{a}^{-1}

                \end{array} \right )
\end{equation}
where $s\in \R$,  $b\in \C$. $N$ is the Heisenberg
group represented  by matrices
$$
   \left ( \begin{array}{ccc}

                        1       &    -2i\bar{b}    &       s-ib\bar{b}  \\

                        0   &    1     &       b \\

                        0      &  0    &       1

                \end{array} \right ).
$$ 

%
\section{The Cartan connection}

The bundle $Y\rightarrow M$ is a $H$-principal bundle and we can interpret  Theorem \ref{structure} by introducing the concept of a Cartan connection.
In fact, one can represent the structure equations of Theorem \ref{structure}
  as a matrix equation whose entries are differential forms. The forms are disposed in the Lie algebra $\mathfrak su(2,1)$ as




$$\pi=\left ( \begin{array}{ccc}
                        -\frac{1}{2}\phi-\frac{1}{3}\omega^1_1     &   -i\phi^{\bar 1}   &   -\frac{1}{4}\psi  \\

                     \omega^1          &   \frac{2}{3}\omega^{1}_{1}    &    \frac{1}{2}\phi^1  \\

                      2\omega         &  2i\omega^{\bar 1}   &      \frac{1}{2}\phi-\frac{1}{3}\omega^1_1
                \end{array} \right )
$$

It is a simple verification to show that 
\begin{equation}\label{dpi}
d\pi+\pi\wedge\pi=\Pi
\end{equation}
where

$$\Pi=\left ( \begin{array}{ccc}
                        0 &   -i\Phi^{\bar 1}      &  -\frac{1}{4}\Psi   \\

                         0  &  0    &  \frac{1}{2}\Phi^1  \\

                           0  &   0  &    0
               \end{array} \right )
$$

Recall that $X^*(y)=\frac{d }{dt}_{_{t=0}}{y e^{tX}}$ where $e^{tX}$ is the one parameter group generated
by the $X$. 

\begin{dfn}

A Cartan connection on $Y$ is a 1-form $\pi: TY \rightarrow su(2,1)$ satisfying:\newline
1. $\pi_p:T_pY\rightarrow su(2,1)$ is an isomorphism\newline
	2. If $X\in \mathfrak h$ and $X^*\in TY$ is the vertical vector field  canonically associated to $X$ then
	$\pi(X^*)=X.$\newline
	3. If $h\in H$ then  $(R_h)^* \pi=Ad_{h^{-1}}\pi$
\end{dfn}

\begin{thm}[\cite{C,CM}]
The form $\pi$ is a Cartan connection.
\end{thm}

\section{Formulae}\label{for}
In the following we will need the value  of $B=Ad_{h^{-1}}\pi$, where $h\in H$ is as in (\ref{h}). We get
$$B_{31}=2(a\bar a\omega)$$
$$B_{21}=a^2{\bar a}^{-1}\omega^1-2a^2b\omega$$
$$B_{11}=-\frac{1}{3}\omega^1_1-\frac{1}{2}\phi+2ia\bar b\omega^1-2a\bar a(s+ib\bar b)\omega$$
$$B_{22}=\frac{2}{3}(\omega^1_1-3ia\bar b\omega^1-3i\bar a b \omega^{\bar 1}+6ia\bar a b \bar b \omega)$$
$$B_{23}=\frac{1}{2}\left(a\bar a^{-2}\phi^1 -ba\bar a^{-1}\phi+2a\bar a^{-1}b\omega^1_1+2a^2\bar a^{-1}(s-ib\bar b)\omega^1-4iab^2\omega^{\bar 1}-4a^2b(s-ib\bar b)\omega\right)$$
$$B_{13}=-\frac{1}{4}\left( (a\bar a)^{-1}\psi+4ia^{-1}b\phi^{\bar 1} -4i\bar a^{-1}\bar b\phi^1-8ib\bar b\omega^1_1\right.
$$
$$
\left. -8ia\bar b(s-ib\bar b)\omega^1  +8i\bar a  b(s+ib\bar b)\omega^{\bar 1}+4s\phi+8a\bar a(s^2+(b\bar b)^2)\omega \right)$$
Observe that 
$$B_{11}+\frac{1}{2}B_{22}=-\frac{1}{2}(\phi-2ia\bar b\omega^1+2i\bar ab\omega^{\bar 1}+4a\bar as\omega).$$

Also, if $h$ is a section (that is a function $h: M\rightarrow H$), then
\begin{equation}\label{dh}
  h^{-1}dh
  = \left ( \begin{array}{ccc}

                        a^{-1}da       &    -2i(a^{-2}\bar a\bar b da+a^{-1}\bar a d\bar{b})    &      ds+i\bar b db-ibd\bar b+a^{-1}(s-ib\bar b)da+\bar a^{-1}(s+ib\bar b)d\bar a 
                          \\

                        0   &   -a^{-1}da+\bar a^{-1}d\bar a     &       a\bar a^{-1}db+a\bar a^{-2}bd\bar a \\

                        0      &  0    &       -\bar{a}^{-1}d\bar a

                \end{array} \right )
\end{equation}
If 
$$\tilde\pi=R_{h}^*\pi=h^{-1}dh+Ad_{h^{-1}}\pi,$$
and writing
$$
\tilde \pi=\left ( \begin{array}{ccc}
                        -\frac{1}{2}\tilde\phi-\frac{1}{3}\tilde\omega^1_1     &   -i\tilde\phi^{\bar 1}   &   -\frac{1}{4}\tilde\psi  \\

                     \tilde\omega^1          &   \frac{2}{3}\tilde\omega^{1}_{1}    &    \frac{1}{2}\tilde\phi^1  \\

                      2\tilde\omega         &  2i\tilde\omega^{\bar 1}   &      \frac{1}{2}\tilde\phi-\frac{1}{3}\tilde\omega^1_1
                \end{array} \right )
$$ 
we obtain from above
$$\tilde\omega=a\bar a\omega$$
$$\tilde\omega^1=a^2{\bar a}^{-1}\omega^1-2a^2b\omega$$
$$\tilde\phi=\phi-2ia\bar b\omega^1+2i\bar ab\omega^{\bar 1}+4a\bar as\omega+(-a^{-1}da-\bar a^{-1}d\bar a)$$
$$\tilde\omega^1_1=\omega^1_1-3ia\bar b\omega^1-3i\bar a b \omega^{\bar 1}+6ia\bar a b \bar b \omega+\frac{3}{2}( -a^{-1}da+\bar a^{-1}d\bar a)$$
$$\tilde\phi^1=a\bar a^{-2}\phi^1 -ba\bar a^{-1}\phi+2a\bar a^{-1}b\omega^1_1+2a^2\bar a^{-1}(s-ib\bar b)\omega^1-4iab^2\omega^{\bar 1}-4a^2b(s-ib\bar b)\omega+2(a\bar a^{-1}db+a\bar a^{-2}bd\bar a )$$
$$\tilde\psi=(a\bar a)^{-1}\psi+4ia^{-1}b\phi^{\bar 1} -4i\bar a^{-1}\bar b\phi^1-8ib\bar b\omega^1_1
 -8ia\bar b(s-ib\bar b)\omega^1  +8i\bar a  b(s+ib\bar b)\omega^{\bar 1}+4s\phi+8a\bar a(s^2+(b\bar b)^2)\omega$$
 $$ -4(ds+i\bar b db-ibd\bar b+a^{-1}(s-ib\bar b)da+\bar a^{-1}(s+ib\bar b)d\bar a ).$$
 We take a particular section $h$ where 
 $$
 b=\frac {\bar a}{a},
 $$
and  write $a=rv$ with $r>0$ and $|v|=1$. The section  $h$ depends on  functions $r$, $v$ and $s$. As $\frac{da}{a}=\frac{dr}{r}+\frac{dv}{v}$, our above formulas can be  writen  as
$$\tilde\omega=r^2\omega$$
$$\tilde\omega^1=rv^3\omega^1-2r^2\omega$$
$$\tilde\phi=\phi-2irv^3\omega^1+2ir\bar v^3\omega^{\bar 1}+4r^2s\omega-2\frac{dr}{r}$$
\begin{equation}\label{simpleh}
\tilde\omega^1_1=\omega^1_1-3ir(v^3\omega^1+{\bar v}^3 \omega^{\bar 1})+6ir^2\omega-3\frac{dv}{v}
\end{equation}
$$
\tilde\phi^1=\frac{1}{r}v^3\phi^1+2\omega^1_1 -\phi+2r(s-i)v^3\omega^1-4ir\bar v^3\omega^{\bar 1}-4r^2(s-i)\omega+2\frac{dr}{r}-6\frac{dv}{v}
$$
$$\tilde\psi=\frac{1}{r^2}\psi-4i\frac 1 rv^3\phi^{ 1} +4i\frac 1 r\bar v^3\phi^{\bar 1}-8i\omega^1_1+4s\phi
 -8rv^3(1+is)\omega^1  -8r\bar v^3(1-is)\omega^{\bar 1}+8r^2(s^2+1)\omega$$
 $$ -4(ds+2s\frac{dr}{r}-6i\frac{dv}{v}).$$
\section{$\R$-Circles}
General references for this section is \cite{C,J}.  The definition (Definition \ref{circles}) of circles by a differential system is due to Cartan, here we introduce a system of differential 
equations describing them (see Theorem \ref{eqcircle}). 
In $S^3$, $\R$-\emph{circles} are the trace of Lagrangian planes of complex hyperbolic space in its boundary $S^3$. To define the analog 
of an $\R$-circle for a CR-structure we begin to impose that the curve is horizontal.  In fact 
we will define the curve in the fiber bundle $Y$ by means of a differential system.  
The projection of an integral curve will be an $\R$-circle.

Observe that a curve $\tilde\gamma$ in $Y$ projects to a horizontal curve $\gamma$ in $M$ if and only if
$$
\omega(\dot{\tilde{\gamma}})=0. 
$$
Indeed,  as $\omega$ is a tautological form of $Y$ over $E$,
$$
\omega(\dot{\tilde{\gamma}})=\lambda\theta(\dot\gamma). 
$$
where  $\lambda$ is positive and $\theta$ is a contact form.

Equation \ref{domega} can be written
$$
d\omega = \omega\wedge\phi + 
\frac i 2(\omega^1-\omega^{\bar 1})\wedge(\omega^1+\omega^{\bar 1}) .
$$ 
In order to obtain an integrable system we add
$$
(\omega^1-\omega^{\bar 1})=0. 
$$
Observe that we could have added instead $(\omega^1-\omega^{\bar 1})=0$ and obtain an equivalent system.
This differential system is still not integrable. It follows from the structure equations in Theorem \ref{structure} that 
$$
d(\omega^1-\omega^{\bar 1})=\frac{1}{2}(\omega^1-\omega^{\bar 1})\wedge\phi+
(\omega^1+\omega^{\bar 1})\wedge\omega^1_1+\omega\wedge(\phi^1-\phi^{\bar 1}),
$$
so as $\omega^1+\omega^{\bar 1}$ is not null, we should add the equation
$$\omega^1_1=0.$$
We rewrite the equation of $d\omega^1_1$ in \ref{structure} as 
$$
d\omega^1_1-\frac{3}{2}i\omega^1\wedge( \phi^{\bar 1}+\phi^{1})-
\frac{3}{2}i(\omega^1-\omega^{\bar 1})\wedge \phi^{ 1}=0\ \  \mod \omega .
$$
In order to obtain an integrable system we add  the equation
$$\phi^1+\phi^{\bar 1}=0.$$
Finally, from the equation of $d\phi^1$ in Theorem \ref{structure} we get $$(\omega^1+\omega^{\bar 1})\wedge \psi=0\ \ \ \mod \omega,\ \omega^1-\omega^{\bar 1}, \ \omega^1_1, \  \phi^1+\phi^{\bar 1}$$
so, as above, in order to obtain an integrable system, we impose
$$\psi=0.$$
Observe that  the structure equations shows that the system is integrable because $d\psi$ is already in the ideal
generated by the previous forms.
We have shown 
\begin{prop}
The differential sistem on $Y$ given by
\begin{equation}\label{eqncircle}
\omega=\omega^1-\omega^{\bar 1}=\omega^1_1=\phi^1+\phi^{\bar 1}=\psi=0
\end{equation}
is integrable
\end{prop}
\begin{dfn}\label{circles}
The curves in the CR-manifold $M$ that are projections of integral manifolds of the differential system \ref{eqncircle} are called $\R$-\emph{circles}.
\end{dfn}

To find the equations of $\R$-circles directly in the CR manifold $M$
 we proceed as in \cite{BS} pg 164 for the case of chains.
Take a section $\sigma : M \rightarrow Y$. Suppose $\gamma:I\rightarrow M$ is a circle, and let's apply a 
transformation  $R_{h(t)}$ on $\sigma\gamma(t)$ such that  $R_{h(t)}\sigma\gamma(t)$ is inside an integral leaf of \ref{eqncircle}. Then
$$
\tilde \pi= (R_{h(t)}\sigma \gamma (t))^*\pi=
h^{-1}(t)h'(t)+Ad_{h^{-1}(t)}\gamma(t)^*(\sigma^*\pi).
$$
Taking account of the formulas \ref{simpleh},  it follows from $\tilde\omega^1-\tilde\omega^{\bar 1}=0$ that $rv^3\omega^1=r\bar v^3\omega^{\bar 1}$, so
\begin{equation}\label{v}
v^3\omega^1=\bar v^3\omega^{\bar 1},
\end{equation}
 and therefore $v$ is determined up to a sixth root of unity.

The following equation $\tilde\omega^1_1=0$ says that
\begin{equation}\label{r}
0=\omega^1_1-3ir(v^3\omega^1+{\bar v}^3 \omega^{\bar 1})-3\frac{dv}{v}
\end{equation}
which determines $r$.

The next equation is 
$$0=\tilde\phi^1+\tilde\phi^{\bar 1}=(\frac{1}{r}v^3\phi^1+2\omega^1_1 -\phi+2r(s-i)v^3\omega^1-4ir\bar v^3\omega^{\bar 1}+2\frac{dr}{r}-6\frac{dv}{v})
$$
$$
+(\frac{1}{r}\bar v^3\phi^{\bar 1}-2\omega^1_1 -\phi+2r(s+i)\bar v^3\omega^{\bar 1}+4ir v^3\omega^{1}+2\frac{dr}{r}+6\frac{dv}{v})
$$
so we obtain
$$
\frac{1}{r}(v^3\phi^1+\bar v^3\phi^{\bar 1})
 -2\phi+2rs(v^3\omega^1+\bar v^3\omega^{\bar 1})+4\frac{dr}{r}=0.
 $$ 
 If we pose
 $$
 z=r^2s
 $$
 we get
\begin{equation}\label{z}
dr=-\frac{1}{4}(v^3\phi^1+\bar v^3\phi^{\bar 1})
 +\frac 1 2r\phi-\frac 1 2z(v^3\omega^1+\bar v^3\omega^{\bar 1})
 \end{equation}
wich determines $z$, and so $s$.

Now the last equation of the integral system is $\tilde\psi=0$, or
$$
0=\frac{1}{r^2}\psi-4i\frac 1 r(v^3\phi^{ 1} -\bar v^3\phi^{\bar 1})+4s\phi
 -8r(v^3\omega^1  +\bar v^3\omega^{\bar 1})
 -4(ds+2s\frac{dr}{r}-6i\frac{dv}{v}),
$$
and from this equation we get
$$
dz=\frac{1}{4}\psi-ir(v^3\phi^{ 1} -\bar v^3\phi^{\bar 1})+z\phi-2r^3(v^3\omega^1  +\bar v^3\omega^{\bar 1})+6ir^2\frac{dv}{v},
$$
If we replace equation \ref{r}  we obtain
\begin{equation}\label{dz}
dz=\frac{1}{4}\psi-ir(v^3\phi^{ 1} -\bar v^3\phi^{\bar 1})+z\phi+2ir^2\omega^1_1+4r^3(v^3\omega^1  +\bar v^3\omega^{\bar 1}) .
\end{equation}
Consider 
$$\gamma'(t)=cZ_1+\bar cZ_{\bar 1},$$ where $Z_1$, $Z_{\bar 1}$, $Z_0$ are dual vectors to $\omega^1$, $\omega^{\bar 1}$, $\omega$. It follows from \ref{v} that $v^3c=\bar v\bar c^3$, so $v^3c=\rho$, with $\rho$ real.
Therefore $c=\rho \bar v^3$, and $(v^3\omega^1  +\bar v^3\omega^{\bar 1})(\gamma')=2\rho$.
If we pose 
$$u=\bar v^3,$$
 we obtain

\begin{thm}\label{eqcircle}If $\gamma :I\rightarrow M$ is an $\R$-circle, and if we write $\gamma'=\rho (u Z_1+\bar u Z_{\bar 1}),$ with $\rho\neq 0$ and $|u|=1$, then $\gamma$ satisfies the following system of differential equations:\newline
$$
\omega=0
$$
$$
\bar u\omega^1-u\omega^{\bar 1}=0
$$
$$
\frac{du}{u}+\omega^1_1-6ir\rho dt=0
$$
$$
dr+\frac{1}{4}(\bar u\phi^1+u\phi^{\bar 1})
 -\frac r 2\phi+z\rho dt=0
$$
$$
dz-\frac{1}{4}\psi+ir(\bar u\phi^{ 1} -u\phi^{\bar 1})-z\phi-2ir^2\omega^1_1-8r^3\rho dt =0.$$ 
\end{thm}

Observe that for each $p\in M$ one can determine a circle passing through $p$ specifying as initial conditions $u$ (with $\vert u\vert =1$), $r>0$ and $z\in \R$. 


\begin{ex}[The quadric, see \cite{J}]
\end{ex}
Consider the \emph{quadric} $Q$ in $\C^2$ defined as the null set of
$$\frac{1}{2i}(z_2-\bar z_2)=z_1\bar z_1.$$
Take as 
$$Z^1=dz_1,$$ 
$$Z^{\bar 1}=d\bar z_1$$ 
and 
$$Z^0=\frac{1}{2}\left(dx_2-i\bar z_1dz_1+iz_1d\bar z_1\right ),$$ 
where $z_2=x_2+iy_2$. The dual fields are $$Z_1=\frac{\partial}{\partial z_1}+i\bar z_1 \frac{\partial}{\partial x_2},$$ 
$$Z_{\bar 1}=\frac{\partial}{\partial \bar z_1}-i z_1 \frac{\partial}{\partial x_2}$$  
and 
$$Z_0=2\frac{\partial}{\partial x_2}.$$
We have
$$dZ^0=iZ^1\wedge Z^{\bar 1}$$ 
and 
$$dZ^1=0.$$ 
Considering this section $\sigma : M\rightarrow Y$ defined by $Z^0,Z^1,Z^{\bar{1}}$, we obtain from Cartan equations  that $\sigma^*\phi=\sigma^*\omega^1_1=\sigma^*\phi^1=\sigma^*\psi=0.$ In particular $Q^1_{\bar 1}=0.$
Replacing this values in equations of Theorem \ref{eqcircle} we obtain
the differential equations for the chains on the quadric, 
\begin{equation}\label{uQ}
\frac{u'}{u}-6i\rho r=0,
\end{equation}
\begin{equation}\label{rQ}
r'+z \rho =0
\end{equation}
and
\begin{equation}\label{sQ}
z'- 8\rho r^3=0.
\end{equation}
One verifies that the curve 
$$\gamma(t)=(x_1(t),y_1(t),x_2(t))$$
given by
$$x_1(t)=\frac{2t (1-t^2)}{1+6 t^2+t^4},\ \ 
y_1(t)=-\frac{(1-t^2) (1+t^2)}{1+6 t^2+t^4},\ \
x_2(t)=-\frac{4t(1+t^2)}{1+6 t^2+t^4}$$
is an $\R$-circle 
corresponding to $z_1(0)=-i$, $x_2(0)=0$, $u(0)=1$, $r(0)=1$ and $z(0)=0$. 
\section{ Circle preserving diffeomorphisms are CR or conjugate CR}

\begin{thm}\label{main}Let $M$ and $\tilde M$ be CR manifolds. If $f:M\rightarrow \tilde M$ is a (local) diffeomorphism such that 
for every $\R$-circle $\gamma$ in $M$, $f\circ \gamma$ is an $\R$-circle in $\tilde M$, then $f$ is a CR diffeomorphism or a conjugate CR diffeomorphism.
\end{thm}

\Pf
Consider the structures $Y$ and $\tilde Y$ on $M$ and $\tilde M$ respectively, and sections $\sigma:M\rightarrow Y$ 
and $\tilde\sigma:\tilde M\rightarrow \tilde Y$. Taking the pull back by $\sigma$ and $\tilde\sigma$ of 
connection forms on $Y$ and $\tilde Y$ respectively, we obtain forms $\omega$, $\omega^1$, $\omega^1_1$, $\phi^1$, $\psi$ on $M$, and forms
 $\tilde\omega$, $\tilde\omega^1$, $\tilde\omega^1_1$, $\tilde\phi^1$,  $\tilde\psi$ on $\tilde M$. Using these forms, we can write circle equations on $M$ and $\tilde M$, as in theorem \ref{eqcircle}. 
We now consider forms $\tilde\omega$, $\tilde\omega^1$, $\tilde\omega^1_1$, $\tilde\phi^1$, 
 $\tilde\psi$ as forms on $M$ (by a slight abuse of notions they denote the forms 
 $f^*\tilde\omega$, $f^*\tilde\omega^1$, $f^*\tilde\omega^1_1$, $f^*\tilde\phi^1$,  $f^*\tilde\psi$). 

As tangent vectors to circles generate $D$, then $\tilde\omega=0$ on $D$, so we obtain $\tilde\omega=\lambda\omega$.  That is, the diffeomorphism is
a contactomorphism. We write
\begin{equation}\label{ab}
\tilde\omega^1= \alpha\omega^1+\beta\omega^{\bar 1}+c\omega
\end{equation}
and
$$
\tilde\omega^{\bar 1}=\bar \beta \omega^1+\bar  \alpha\omega^{\bar 1}+\bar c\omega,
$$
where $\alpha,\beta,c$ are functions on $M$. Observe that as $\tilde\omega^1\wedge\tilde\omega^{\bar 1}|_D=(\alpha\bar \alpha-\beta\bar \beta)\omega^1\wedge\omega^{\bar 1}|_D$ is never null on $D$, we get
$$\alpha\bar \alpha-\beta\bar \beta\neq 0$$
everywhere on $M$.   We have to show that $\beta=0$ (for $f$ to be a CR diffeomorphism) or $\alpha=0$ (for $f$ to be a conjugate CR diffeomorphism).
The idea is to use the arbitrary initial conditions (so vary $u$, $r$ and $z$) to obtain enough information on those functions. 

The computations we are going to do are on a circle $\gamma$, so we can assume from theorem \ref{eqcircle} that 
$$\omega=\tilde\omega=\omega^1-u^2\omega^{\bar 1}=\tilde\omega^1-\tilde u^2\tilde\omega^{\bar 1}=0.$$
Applying this in \ref{ab} we obtain
$$(\alpha u^2+\beta)\omega^{\bar 1}=\alpha\omega^1+\beta\omega^{\bar 1}=\tilde\omega^1=\tilde u^2\tilde\omega^{\bar 1}=\tilde u^2(\bar \beta \omega^1+\bar \alpha\omega^{\bar 1})=\tilde u^2(\bar \beta u^2+\bar \alpha)\omega^{\bar 1},$$
therefore
\begin{equation}\label{tildev}
\tilde u=\left(\frac{\alpha u^2+b}{\bar \beta u^2+\bar \alpha}\right)^{\frac{1}{2}}.
\end{equation}
Define $h:\mathbb R\times S^1\rightarrow \mathbb C$ by 
\begin{equation}\label{hty}
h(t,u)=\left(\frac{\alpha(\gamma(t))u^2+\beta(\gamma(t))}{\bar \beta(\gamma(t)) u^2+\bar \alpha(\gamma(t))}\right)^{\frac{1}{2}}
\end{equation}
for $t\in \mathbb R$ and $u\in S^1\subset \mathbb C$.  
Applying again Theorem \ref{eqcircle} we can write
$$d\tilde u=\tilde u(-\tilde\omega^1_1+6i\tilde r\tilde\rho dt)=h(-\tilde\omega^1_1+6i\tilde\rho\tilde rdt)$$
and taking the derivative of \ref{hty}, we obtain
$$d\tilde u=\frac{\partial h}{\partial u}du+\frac{\partial h}{\partial t}dt=\frac{\partial h}{\partial u}u(-\omega^1_1+6i r\rho dt)+\frac{\partial h}{\partial t}dt.
$$
From the equality of the right sides on the two lines above we obtain
\begin{equation}\label{tildey}
\tilde r=Ar+B
\end{equation}
where
\begin{equation}\label{A}
A(t,u)=\frac{\rho u}{ \tilde\rho h(t,u)}\frac{\partial h}{\partial u}(t,u)
\end{equation}
and $B$ is a function of $t,u$.
Taking the derivative of \ref{tildey} and replacing the circle equation 
$$
d\tilde r=-\frac{1}{4}(\bar{ \tilde u}\tilde \phi^1+\tilde u\tilde \phi^{\bar 1})
 +\frac {\tilde r} 2\tilde \phi-\tilde z\tilde \rho dt
$$
we obtain
$$
-\frac{1}{4}(\bar{ \tilde u}\tilde \phi^1+\tilde u\tilde \phi^{\bar 1})
 +\frac {1} 2(Ar+B)\tilde \phi-\tilde z\tilde \rho dt=
$$
$$=A(-\frac{1}{4}(\bar u\phi^1+u\phi^{\bar 1})
 +\frac r 2\phi-z\rho dt)+(\frac{\partial A}{\partial u}r+\frac{\partial B}{\partial u})(-u\omega^1_1+6i ur\rho dt)+
 (\frac{\partial A}{\partial t}r +\frac{\partial B}{\partial t}) dt,
$$
and  solving for $\tilde z$ we obtain
\begin{equation}\label{tildez}
\tilde z=Cz+Dr^2+Er+F,
\end{equation}
where, in a short form,
\begin{equation}\label{C}
C=\frac{ A\rho}{ \tilde \rho},
\end{equation}
\begin{equation}\label{D}
D=-6i\frac{ \rho u}{ \tilde \rho}\frac{\partial A}{\partial u}.
\end{equation}
Taking the derivative of \ref{tildez} we obtain
$$d\tilde z=C\left(\frac{1}{4}\psi-ir(\bar u\phi^{ 1} -u\phi^{\bar 1})+z\phi+2ir^2\omega^1_1+8r^3\rho dt \right)
$$
$$+(2rD+E)\left(-\frac{1}{4}(\bar u\phi^1+u\phi^{\bar 1})
 +\frac r 2\phi-z\rho dt\right)+(\frac{\partial C}{\partial u}z+\frac{\partial D}{\partial u}r^2+\frac{\partial E}{\partial u}r+\frac{\partial F}{\partial u})(-u\omega^1_1+6i ur\rho dt)
$$
$$
+(\frac{\partial C}{\partial t}z+\frac{\partial D}{\partial t}r^2+\frac{\partial E}{\partial t}r+\frac{\partial F}{\partial t})dt.
$$
By theorem \ref{eqcircle} and \ref{tildey}, \ref{tildez} we have
$$
d\tilde z=\frac{1}{4}\tilde \psi-i(Ar+B)(\bar{\tilde  u}\tilde \phi^{ 1} -\tilde u\tilde \phi^{\bar 1})+(Cz+Dr^2+Er+F)\tilde \phi+2i(Ar+B)^2\tilde \omega^1_1+8(Ar+B)^3\tilde\rho dt.
$$
Replacing this in the above equation, we get a polynomial equation in $z$ and $r$:
\begin{equation}\label{poly}
A_{30}r^3+A_{11}rz+A_{20}r^2+A_{10}r+A_{00}=0
\end{equation}
where
\begin{equation}\label{A11}
A_{11}=( D-3 i  u \frac{\partial C}{\partial u})\rho=-9iu\frac{\rho^2}{\tilde\rho}\frac{\partial A}{\partial u}.
\end{equation}
As $r$ and $z$ can take arbitrary values in $\mathbb R$ we obtain that all coeficients of \ref{poly} are null. If we replace \ref{A}, 
 \ref{C} and \ref{D} in  $A_{11}=0$ we obtain
$$18i\frac{\rho^3}{\tilde\rho^2}\frac{ \left(u^2 \alpha\bar \beta-\bar u^2\bar \alpha
   \beta\right) (\alpha \bar \alpha-\beta
   \bar \beta)}{\left(\alpha u+\beta\bar u\right)^2
   \left(\bar \beta u+\bar \alpha\bar u\right)^2}=0.
   $$
As $u$ is any complex number such that $|u|=1$ and $\alpha \bar \alpha-\beta
   \bar \beta\neq 0$, we obtain $\alpha\bar \beta=0$. It follows from \ref{ab} that if $\beta=0$ then $f$ is a CR diffeomorphism and if $\alpha=0$ we get that $f$ is a conjugate CR diffeomorphism.
   
   \EPf

\end{document}